\documentclass[11pt,dvipsnames]{amsart}
\usepackage{xy}
\usepackage{graphicx, amsmath, amssymb, amsthm,amsfonts}
\usepackage[mathscr]{eucal}
\usepackage[pagebackref, colorlinks, citecolor=Red, linkcolor=NavyBlue, urlcolor=NavyBlue]{hyperref}
\hypersetup{linktoc=all}
\usepackage{comment}
\usepackage{hyperref}
\usepackage[capitalise, nameinlink]{cleveref}
\usepackage{microtype}
\usepackage{tikz-cd}
\usepackage{spectralsequences}

\newcommand{\Z}{\mathbb{Z}}
\newcommand{\A}{\mathcal{A}}
\newcommand{\F}{\mathbb{F}}
\newcommand{\p}{\mathfrak{p}}
\newcommand{\q}{\mathfrak{q}}

\crefname{lemma}{Lemma}{Lemmas}
\crefname{theorem}{Theorem}{Theorems}
\crefname{definition}{Definition}{Definitions}
\crefname{proposition}{Proposition}{Propositions}
\crefname{remark}{Remark}{Remarks}
\crefname{corollary}{Corollary}{Corollaries}
\crefname{equation}{Equation}{Equations}
\crefname{construction}{Construction}{Constructions}
\crefname{ex}{Example}{Examples}
\crefname{example}{Example}{Example}
\crefname{appsec}{Appendix}{Appendices}
\crefname{subsection}{Subsection}{Subsections}
\crefname{letterthm}{Theorem}{Theorems}

\newtheorem{theorem}{Theorem}[section]
\newtheorem{letterthm}{Theorem}

\newtheorem{lettercor}[letterthm]{Corollary}

\newtheorem{lemma}[theorem]{Lemma}

\newtheorem{proposition}[theorem]{Proposition}
\newtheorem{corollary}[theorem]{Corollary}

\newtheoremstyle{BoldRemark} 
{10pt}                    
{10pt}                    
{\upshape}                   
{}                           
{\bfseries}                  
{.}                          
{.5em}                       
{}  
\theoremstyle{BoldRemark}
\newtheorem{remark}[theorem]{Remark}
\newtheorem{example}[theorem]{Example}
\newtheorem{definition}[theorem]{Definition}
\newtheorem{letterdef}[letterthm]{Definition}

\newtheorem{question}[theorem]{Question}

\usepackage[textwidth=25mm, textsize=tiny]{todonotes}

\newcommand{\Res}{\mathrm{Res}}
\newcommand{\res}{\mathrm{res}}
\newcommand{\Spec}{\mathrm{Spec}}

\newcommand{\Nm}{\mathrm{Nm}}

\newcommand{\Coind}{\mathrm{CoInd}}

\newcommand{\Tr}{\mathrm{Tr}}

\newcommand{\FP}{\mathrm{FP}}

\renewcommand{\subset}{\subseteq}
\newcommand{\ghost}{\Gamma}

\newcommand{\Mod}{\mathrm{Mod}}

\newcommand{\Tamb}{\mathrm{Tamb}}

\usepackage[margin=1in]{geometry}

\title{The subgroup stratification of Nakaoka spectra}

\author[N. Wisdom]{Noah Wisdom}
\address{Department of Mathematics,
         Northwestern University,
         Evanston, IL, 60208
         USA}
\email{NoahAnkney2026@u.northwestern.edu}



\date{}

\begin{document}
\begin{abstract}
    We construct a stratification on the Nakaoka spectrum of any $G$-Tambara functor indexed by the poset of subgroups of $G$. When $G$ is Dedekind, we show that the $H$th stratum of the Nakaoka spectrum of the Burnside $G$-Tambara functor is closed and not open; this provides examples of \'{e}tale maps which induce closed, non-open maps on Nakaoka spectra. By computing the strata on the ghost of a $C_p$-Tambara functor we obtain many more examples of \'{e}tale maps of Tambara functors for which the induced map on Nakaoka spectra is not open, in contrast to the non-equivariant world. We also compute the strata of all fixed-point Tambara functors.
\end{abstract}

\maketitle
\tableofcontents

\section{Introduction}

Tambara functors are equivariant analogues of rings appearing in equivariant homotopy theory \cite[Section 7.2]{Bru07}, \cite[Theorem 1.4]{AB18}, \cite{CHLL24,GKM26}, and in many other purely algebraic situations \cite[Section 3]{Tam93}, \cite{Bru05,Nak10,Sul23}. They come equipped with a Nakaoka spectrum analogous to the Zariski spectrum, the collection of prime ideals with a certain topology. One of the basic problems in equivariant algebra is to compute the Nakaoka spectrum of a Tambara functor. For example, the computation of the Nakaoka spectrum of the Burnside $G$-Tambara functor $\A_G$---the analogue of the integers in equivariant algebra---took over a decade and spanned the work of many authors \cite{Nak11a,Nak14,CG23,4DS,4DSC}. Moreover, Tambara functor algebraic geometry is in its infancy; there is not yet even an accepted definition of the structure sheaf of a Nakaoka spectrum. The present article studies the topology on Nakaoka spectra, hopeful that it will shed light on how to define a useful analogue of the structure sheaf in the Tambara setting.

Given a subgroup $H$ of $G$, there are various change-of-group functors relating the category of $H$-Tambara functors and the category of $G$-Tambara functors: of particular importance to us is the restriction 
\[ 
    \Res_H^G : G \text{-} \Tamb \rightarrow H \text{-} \Tamb 
\] 
which is left adjoint to coinduction $\Coind_H^G$. Since $\Res_H^G$ preserves initial objects, we have $\Res_H^G \A_G \cong \A_H$. Previous work \cite{SSW24,Wis25a} suggests that $\Coind_H^G$ has extremely nice commutative-algebraic properties. We extend this philosophy to ideals:

\begin{letterthm}(\cref{thm:coinduction-induces-bijection-on-Tambara-ideals})
    $\Coind_H^G$ induces an order-preserving bijection between Tambara ideals of $R$ and Tambara ideals of $\Coind_H^G R$ which restricts to a bijection between prime ideals. Consequently we obtain a natural homeomorphism
    \[ 
        \Spec(\Coind_H^G R) \cong \Spec(R) 
    \] 
    of Nakaoka spectra for all $H$-Tambara functors $R$.
\end{letterthm}

The adjunction unit 
\[ 
    R \rightarrow \Coind_H^G \Res_H^G R 
\] 
thus induces compatible maps of Nakaoka spectra.

\begin{letterdef}
    Let $R$ be a $G$-Tambara functor. We define the \emph{subgroup stratification} of $\Spec(R)$ as the stratification indexed by the poset of subgroups of $G$ whose $H$th stratum is the image of a certain continuous map
    \[ 
        \Spec(\Res_H^G R) \rightarrow \Spec(R) .
    \] 
\end{letterdef}

\begin{remark}
    One of the motivations of \cite{4DSC} is to determine the target of a conjectured comparison map from the Balmer spectrum of a $G$-tensor-triangulated category to the Nakaoka spectrum of the endomorphism Tambara functor associated to the unit object. Assuming such a comparison map exists, we obtain a stratification on the Balmer spectrum by taking the preimage of the subgroup stratification under the comparison map. This stratification on Balmer spectra of $G$-tensor-triangulated categories should then enjoy the same naturality properties as the subgroup stratification on Nakaoka spectra.
\end{remark}

Recall that a Dedekind group is a group for which all subgroups are normal.

\begin{letterthm}(\cref{thm:strata-of-Burnside-G-Dedekind,prop:stratum-of-ghost,prop:subgroup-strata-for-MRC})
    We have the following explicit descriptions of strata: 
    \begin{enumerate}
        \item For the Burnside functor $\A_G$, if $G$ is Dedekind, the $H$th stratum of $\Spec(\A_G)$ is the basic closed subset $V(\p_{H,0})$ which is not open.
        \item If $R$ is a $C_p$-Tambara functor and $\ghost(R)$ its ghost \cite[Definition 7.3]{4DS}, the $e$th stratum of $\Spec(\ghost(R))$ is the subset of points $(\p;\Phi^{C_p}(R))$ where $\p \in \Spec(R(C_p/e)^{C_p})$.
        \item For $G$ arbitrary, if $Q$ is a $G$-ring, then the $H$th stratum of $\Spec(\FP(Q))$ is the subspace of prime ideals $\p$ of $\FP(Q)$ such that $\p(G/H)$ is a prime ideal of $Q^H$. Equivalently, it is the subspace of those $\p$ for which $\p(G/e)$ is $H$-prime.
    \end{enumerate}
\end{letterthm}

\begin{figure}[h]
  \centering
  \includegraphics[width=\textwidth]{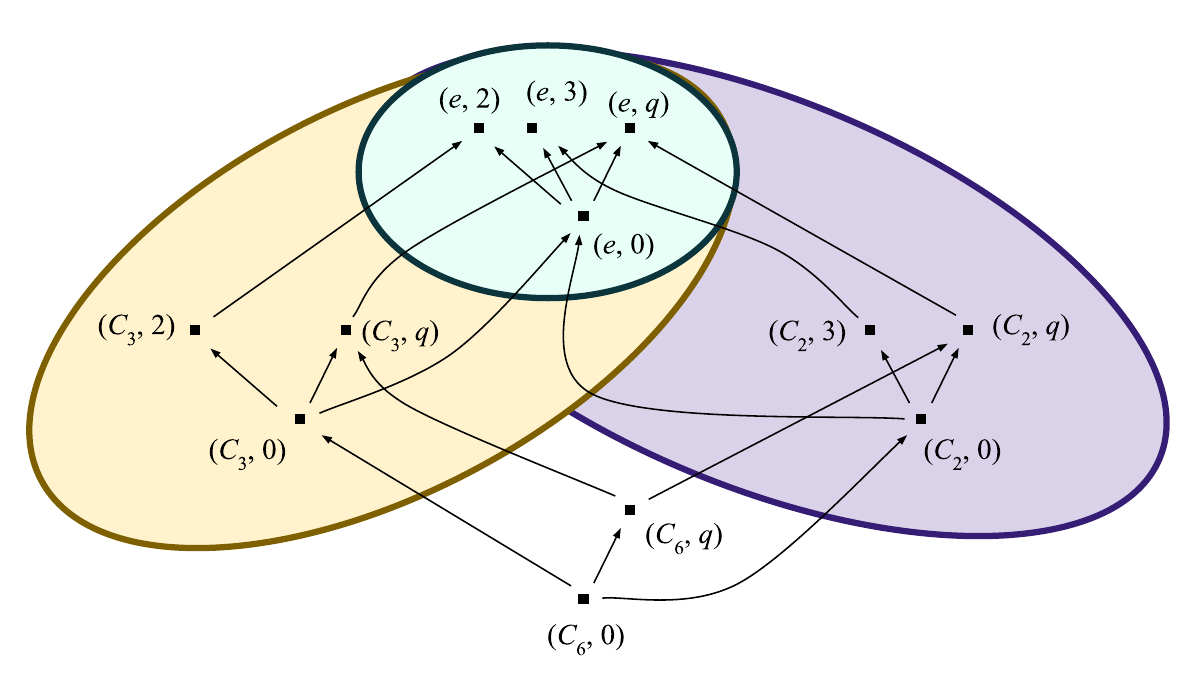}
  \caption{\textbf{The subgroup stratification on $\Spec(\A_{C_6})$.} A point labelled $(H,p)$ represents the prime ideal $\p_{H,p}$, and there is a path of arrows from point $\p$ to point $\p'$ precisely when there is a proper inclusion $\p \subsetneq \p'$. $q \neq 2, 3$ denotes a generic prime. The oval at the top is the $e$th stratum. The left oval is the $C_3$th stratum; note that it contains the $e$th stratum and is homeomorphic to $\Spec(\A_{C_3})$. Similarly, the right oval is the $C_2$th stratum. The entire space is the $C_6$th stratum.}
  \label{fig:subgroup_strat_on_Spec_A_C_6}
\end{figure}

For example, \cref{fig:subgroup_strat_on_Spec_A_C_6} shows the subgroup stratification on the Nakaoka spectrum of the Burnside Tambara functor when $G = C_6$.

\begin{remark}
    When $G = C_p$, by \cref{thm:strata-of-Burnside-G-Dedekind} the $e$th stratum of $\Spec(A_{C_p})$ is homeomorphic to $\Spec(\Z)$, and its complement is given as a set by $\{ \p_{C_p,q} | q \in \Spec(\Z), q \neq p \}$. In particular, one of the most interesting parts of Nakaoka's description of $\Spec(\A_{C_p})$ (\cite[Theorem 1.1]{Nak14}, reviewed in \cref{subsec:Burnside_example}) is the collision $\p_{e,p} = \p_{C_p,p}$, and this manifests in the subgroup stratification perspective as a kind of trans-stratum phenomenon: the ideal determined by the notation $\p_{C_p,p}$ is unexpectedly in the $e$th stratum.

    We note that this collision in $\Spec(\A_{C_p})$ is also related to the transchromatic blueshift phenomenon appearing in the Balmer spectrum of $G$-spectra \cite{BS17,BHNNNS,BGH20}. One might also wonder about a connection to the topology on the Balmer spectrum studied in \cite{PSW22}.
\end{remark}

From \cite[Theorem A]{Wis25c} $\A_G \rightarrow \Coind_H^G \A_H$ is \'{e}tale, i.e. it is finitely presented, flat, and the genuine K\"{a}hler differentials of \cite[Definition 4.1]{Hil17} vanish. Thus from \cref{thm:strata-of-Burnside-G-Dedekind} we immediately obtain the following.

\begin{lettercor}(\cref{cor:ex-of-etale-map-which-is-not-open})
    Let $G$ be an arbitrary finite group. There exists an \'{e}tale map of $G$-Tambara functors for which the induced map on Nakaoka spectra is not open. If $G$ is Dedekind, then the induced map on Nakaoka spectra is closed.
\end{lettercor}

When $G = C_p$ we construct many more examples of this form using the ghost construction $\Gamma(R)$ of a $C_p$-Tambara functor $R$ (\cite[Definition 7.3]{4DS} reviewed in \cref{subsec:ghost_example}).

\begin{letterthm}(\cref{cor:ghost-strata-are-usually-not-open})
    Let $R$ be a $C_p$-Tambara functor such that the transfer is not surjective. Then the canonical map
    \[ 
        \ghost(R) \rightarrow \Coind_e^{C_p} R(C_p/e) 
    \] 
    is \'{e}tale, but the induced map on Nakaoka spectra is not open.
\end{letterthm}

Finally, we give novel necessary conditions for a Tambara functor to be domain-like (\cref{def:domain-like}).

\begin{letterthm}[\cref{thm:strata-usually-not-open-in-domain-like}]
    Let $G$ be an arbitrary group and $R$ any domain-like Tambara functor such that the restriction $R(G/G) \rightarrow R(G/H)$ is not injective. Then the $H$th stratum of $\Spec(R)$ is not open.
\end{letterthm}

\begin{remark}
    Let $R$ be a $C_p$-Tambara functor. Then we may consider the classical Nakaoka spectrum $\Spec(R)$. On the other hand, the $C_p$-action on $R(C_p) = \Res_e^{C_p} R$ determines a $C_p$-action on the Zariski spectrum $\Spec(\Res_e^{C_p} R)$. Motivated by Elmendorf's theorem, which does not apply literally here, one is led to ask for the existence of an algebraically defined $C_p$-space closely related to $\Spec(\Res_e^{C_p} R)$ (or possibly a family of deformations thereof) whose fixed points are the Nakaoka spectrum of $R$. Such a space could support well-behaved Tambara functor algebraic geometry.

    This question is closely related to the following question, which the author learned of from Mike Hill: how can one recover the Nakaoka spectrum of $R$ from the geometric information of the affine schemes $\Spec(R(G/H))$?
\end{remark}

One might also attempt to define a topology on the set of prime ideals of a Tambara functor by declaring all \'{e}tale maps to be open. This perspective suggests a Tambara analogue of the classical notion of algebraic spaces as developed by Artin.

\vspace{3mm} \textbf{Acknowledgments.} The author thanks Mike Hill for the encouragement to spend more time thinking about the subgroup stratification, and for many other insightful suggestions. Additionally, the author thanks David Chan and David Mehrle for helpful conversations. Finally, the author thanks the anonymous referee for numerous suggestions which greatly improved the quality of the article, including the suggestion to compare to $\Spec(A(G))$ in \cref{rem:cf_4DSC_Theorem_4.19}.

\section{Recollections on Nakaoka spectra and coinduction}

Any commutative ring possesses a Zariski spectrum, the space of prime ideals. Similarly, a Tambara functor gives rise to a Nakaoka spectrum, a space of prime ideals. We follow \cite{Nak11a}, although see \cite{4DS} for additional excellent exposition.

\begin{definition}[{cf. \cite[Definition 2.1]{Nak11a}}]
	A \emph{Tambara ideal} of a $G$-Tambara functor $R$ is a collection of ideals $I(G/H) \subset R(G/H)$ subject to the conditions
	\begin{enumerate}
		\item $\mathrm{Tr}_K^H \left( I(G/K) \right)  \subset I(G/H)$
		\item $\mathrm{Nm}_K^H \left( I(G/K) \right) \subset I(G/H)$
		\item $\mathrm{Res}_K^H \left( I(G/H) \right) \subset I(G/K)$
		\item $c_g \left( I(G/H) \right) \subset I(G/g^{-1}Hg)$
	\end{enumerate}
	where $c_g$ denotes restriction along right multiplication by $g$, $xH \mapsto xgg^{-1}Hg$.
\end{definition}

Natalie Stewart and David Chan have pointed out to the author that the condition 
\[
	\mathrm{Nm}_K^H \left( I(G/H) \right) \subset \mathrm{Nm}_K^H(0) + I(G/K)
\] 
of \cite[Definition 2.1]{Nak11a} is superfluous when we take the convention that Tambara functors take values only on transitive $G$-sets, rather than all finite $G$-sets, as the norm of $0 \in R(G/e)$ is always zero.

If $R$ is merely a Green functor and we relax the requirement that $I$ be preserved by norms, we obtain the notion of a Green ideal. A Green ideal of a Green functor $R$ is just an $R$-submodule of $R$.

Nakaoka proves that the intersection of Tambara ideals is again a Tambara ideal. Thus given any collection of elements in a Tambara functor, there is a unique smallest Tambara ideal containing them. In particular, we may speak of principal Tambara ideals, and more generally Tambara ideals generated by a set.

\begin{definition}
	The product $I \cdot J$ of two ideals $I$ and $J$ of a Tambara functor is defined as the Tambara ideal generated by the set of elements $I(G/H) \cdot J(G/H)$ in each level $G/H$. A Tambara ideal $\p$ is called prime if $I \cdot J \subset \p$ implies either $I \subset \p$ or $J \subset \p$.
\end{definition}

\begin{example}
    Since the Burnside Tambara functor is the initial $G$-Tambara functor, there is a unique Tambara functor map $\A_G \rightarrow \underline{\Z}$, where $\underline{\Z}$ is the fixed-point $G$-Tambara functor associated to the trivial $G$-action on $\Z$. The composition 
    \[ 
        \A_G(G/H) \rightarrow \underline{\Z}(G/H) = \Z \cong \A_G(G/G) 
    \] 
    is equal to the restriction $\res_H^G$ in $\A_G$; equivalently, it is induced by the map sending a finite $H$-set to its cardinality, using the identification of $\A_G(G/H)$ with the Burnside ring for the group $H$. By \cite[Proposition 2.10]{Nak11a} the collection of kernels 
    \[ 
        \mathrm{Ker}(\A_G(G/H) \rightarrow \Z) 
    \] 
    form a Tambara ideal of $\A_G$, which one may check is prime. In fact, this is the prime ideal $\p_{e,0}$ in the notation of \cite{4DSC}.
\end{example}

Now we may provide Nakaoka's definition of the spectrum of a Tambara functor.

\begin{definition}\cite[Section 4.2]{Nak11a}
    The \emph{Nakaoka spectrum} of a Tambara functor $R$ is the set of prime Tambara ideals of $R$ equipped with the topology generated by closed subsets
    \[ 
        V(I) := \{ \p | I \subset \p \} 
    \] 
    where $I$ runs through all Tambara ideals of $R$.
\end{definition}

It is important to note that the space $\mathrm{Spec}(R)$ is determined entirely by the poset of Tambara ideals of $R$ along with its sub-poset of prime ideals.

We now turn our attention to the functor $\Coind_H^G$ which will play a central role.

\begin{definition}
    Let $R$ be an $H$-Tambara (resp. Green, Mackey) functor. Then we define a $G$-Tambara (resp. Green, Mackey) functor $\Coind_H^G R$ as follows:
    \[ 
        \Coind_H^G R(G/H) := R(\Res_H^G G/H) 
    \] 
    where we employ the convention $R(X \sqcup Y) := R(X) \times R(Y)$. If $Q$ is an $H$-ring, define a $G$-ring by
    \[ 
        \Coind_H^G Q := \mathrm{Fun}(G/H,Q) 
    \] 
    where addition and multiplication are pointwise.
\end{definition}

We will heavily use the following result from \cite{Wis25a}, which we recall here. 

\begin{theorem}[{\cite[Theorem F and Corollary G]{Wis25a}}]
    Let $R$ be an $H$-Tambara functor. $\Coind_H^G$ induces equivalences of categories
    \begin{align*}
        R \text{-} \Mod & \cong ( \Coind_H^G R ) \text{-} \Mod \\
        R \text{-} \mathrm{Alg} & \cong ( \Coind_H^G R ) \text{-} \mathrm{Alg} .
    \end{align*} 
\end{theorem}

Since we will also use a structure theorem for Tambara functors from \cite{Wis25a}, we recall it as well.

\begin{definition}[{\cite[Definition C]{Wis25a}}]\label{def:clarified}
	A $G$-ring $Q$ is \emph{clarified} if it does not contain any idempotents $d$ with both of the following properties: 
	\begin{enumerate}
		\item the isotropy group of $d$ is a proper subgroup of $G$, and
		\item the $G$-orbits of $d$ which are not equal to $d$ are orthogonal to $d$, that is, $gd \neq d$ implies $d \cdot gd = 0$.
	\end{enumerate}
    We call a Tambara functor $R$ clarified if the $G$-ring $R(G/e)$ is.
\end{definition}

\begin{theorem}[{\cite[Theorem D]{Wis25a}}]
	Let $R$ be a Tambara functor such that $R(G/e)$ is Noetherian. Then there exists an isomorphism 
    \[ 
        R \cong \prod_{H \subset G} \mathrm{Coind}_H^G R_H 
    \] 
    such that each $R_H$ is a clarified $H$-Tambara functor.
\end{theorem}

Finally, we will take advantage of a family of examples of \'{e}tale maps of Tambara functors constructed in \cite{Wis25c} in order to record a pathology of \'{e}taleness in Tambara functor commutative algebra.

\begin{theorem}[{\cite[Theorem A]{Wis25c}}]
    Let $R$ be a $G$-Tambara functor. The adjunction unit 
    \[ 
        R \rightarrow \Coind_H^G \Res_H^G R 
    \] 
    is \'{e}tale in the sense that it is flat, finitely presented, and the genuine K\"{a}hler differentials of \cite[Definition 4.1]{Hil17} vanish.
\end{theorem}

\section{The subgroup stratification}

We start by establishing some results on the interaction between products, coinductions, and ideals.

\begin{lemma}\label{lem:product-induces-bijection-on-Tambara-ideals}
	Let $R \cong S \times T$ be an isomorphism of Tambara functors. The projections $R \rightarrow S$ and $R \rightarrow T$ determine a bijection between the set of prime Tambara ideals of $R$ and the disjoint union of the respective sets of prime Tambara ideals of $S$ and $T$.
\end{lemma}

\begin{proof}
    Let $I$ be a prime ideal of $R$, and let $e_S = (1,0)$, $e_T = (0,1)$ in $R(G/G)$. We have 
    \[ 
        (e_S) \cdot (e_T) = (0) \subset I 
    \] 
    so that either $e_S \in I$ or $e_T \in I$. Equivalently, either $I$ is $S \times J$ for some uniquely determined ideal $J$ of $T$, or else $I$ is $J \times T$ for some $J \subset S$. Without loss of generality we are in the former case; since $T/J \cong R/I$ we deduce that $J$ is prime (as primeness of an ideal is a property of the quotient by that ideal). It follows that 
    \[ 
        \Spec(S) \sqcup \Spec(T) \rightarrow \Spec(R) 
    \] 
    is surjective. It is injective because $J$ is uniquely determined by $I$.
\end{proof}

\begin{theorem}\label{thm:coinduction-induces-bijection-on-Tambara-ideals}
	Coinduction induces a bijection between the set of Tambara ideals of $R$ and the set of Tambara ideals of $\Coind_H^G R$. This bijection respects containment and products of ideals, hence restricts to a bijection between prime ideals.
\end{theorem}

\begin{proof}
	To start with, observe that by \cite[Theorem F]{Wis25a} $\Coind_H^G$ induces a bijection between Green ideals of the Green functor underlying $R$ and the Green ideals of the Green functor underlying $\Coind_H^G R$. We will repeatedly use the following description of the inverse of this bijection. First, \cite[Proposition 5.3]{TW95} (which is stated for Mackey functors, but is true by an identical proof for Tambara functors and for modules over a fixed Tambara functor) supplies the equation 
    \[ 
	   \mathrm{Res}_K^G \mathrm{Coind}_H^G M \cong \prod_{g \in K \backslash G / H} \mathrm{Coind}_{K \cap {}^g H}^K \mathrm{Res}_{K \cap {}^g H}^{{}^g H} {}^g M 
    \] 
    when $M$ is a Mackey functor, Tambara functor, or module over a fixed Tambara functor. The proof of \cite[Theorem F]{Wis25a} shows that the inverse to coinduction
    \[ 
        (\Coind_H^G)^{-1} : (\Coind_H^G R) \text{-} \Mod \xrightarrow{\cong} R \text{-} \Mod 
    \] 
    is given by restriction, followed by projection onto the identity double coset factor.
	
	Suppose $I$ is a Green ideal of $R$. We must show that $\Coind_H^G I$ is closed under norms if and only if $I$ is. Suppose $\Coind_H^G I$ is a Tambara ideal. Since $\mathrm{Res}_H^G$ preserves Tambara ideals, we see that $\mathrm{Res}_H^G \Coind_H^G I$ is a Tambara ideal of $\mathrm{Res}_H^G \Coind_H^G R$. Now the latter is a product of $H$-Tambara functors, and one factor is $R$. We have already seen that every Tambara ideal of a product of Tambara functors is a set-theoretic product of Tambara ideals, and the factor of $\mathrm{Res}_H^G \Coind_H^G R$ corresponding to $R$ is $I$. Hence $I$ is a Tambara ideal. On the other hand, if $I$ is closed under norms, then since norms in $\Coind_H^G R$ are given levelwise as the set-theoretic product of ring-theoretic products of norms in $R$, and Green ideals are closed under products, we see that $\Coind_H^G I$ is closed under norms.
	
	Next, we show that this bijection respects products of ideals. Let $I$ and $J$ be Tambara ideals of $R$. Recall that $\Coind_H^G R(G/K)$ is a product of copies of rings $R(H/L)$, and in particular the coinduction of $I$ is given in level $G/K$ by the very same set-theoretic product of ideals $I(H/L) \subset R(H/L)$. The factor of $( \Coind_H^G ( I \cdot J ))(G/K)$ which lives in $R(H/L)$ contains $I(H/L) \cdot J(H/L)$. Thus 
    \[ 
        \Coind_H^G I \cdot \Coind_H^G J \subset \Coind_H^G ( I \cdot J ) 
    \]
    On the other hand, $\Coind_H^G I \cdot \Coind_H^G J$ is generated as a Tambara ideal by the set-theoretic product of ideal products $I(H/L) \cdot J(H/L) \subset R(H/L)$ over the factors of $\Coind_H^G R(G/K)$. The projection of the restriction of $\Coind_H^G I \cdot \Coind_H^G J$ to the identity double coset factor therefore contains $I \cdot J$. Thus 
    \[
        \Coind_H^G ( I \cdot J ) \subset \Coind_H^G I \cdot \Coind_H^G J 
    \] 
    so that in fact equality holds.

	Finally, we must show that $\p$ is prime if and only if $\Coind_H^G \p$ is. This follows from the fact that coinduction induces a bijection between the sets of Tambara ideals which respects both inclusions and ideal products, though we give one direction for clarity. Suppose $\Coind_H^G \p$ is prime and let $I \cdot J \subset \p$. Then since coinduction preserves ideal products and inclusions, we observe $\Coind_H^G I \cdot \Coind_H^G J \subset \Coind_H^G \p$, so that (without loss of generality) $\Coind_H^G I \subset \Coind_H^G \p$. Restricting and projecting on the identity double coset factor preserves inclusion, so $I \subset \p$.
\end{proof}

Since the topology on the Nakaoka spectrum is completely determined by the poset of Tambara ideals along with the sub-poset of prime ideals, we immediately obtain the following.

\begin{corollary}\label{cor:coind-induces-Nakaoka-spec-homeo}
    For any $H$-Tambara functor $R$ we have a natural homeomorphism 
    \[ 
        \Spec(\Coind_H^G R) \cong \Spec(R) 
    \] 
    of Nakaoka spectra.
\end{corollary}

Let $Q$ be a Noetherian ring. There is a one-to-one correspondence between connected components of $\Spec(Q)$ and ``minimal" idempotents in $Q$. If $R$ is a Tambara functor, one may similarly use idempotents in $R(G/e)$ to decompose $\Spec(R)$ as a disjoint union of spaces. One incarnation of this idea is the following.

\begin{theorem}\label{thm:Nakaoka-spectrum-computation}
	Any Nakaoka spectrum $\mathrm{Spec}(R)$ where $R(G/e)$ is Noetherian decomposes as a disjoint union over subgroups $H$ of $G$ of Nakaoka spectra of clarified $H$-Tambara functors.
\end{theorem}

\begin{proof}
    Combine \cite[Theorem D]{Wis25a} with \cref{lem:product-induces-bijection-on-Tambara-ideals,thm:coinduction-induces-bijection-on-Tambara-ideals} .
\end{proof}

We arrive now at our main definition.

\begin{definition}
    Let $R$ be a $G$-Tambara functor. We define the \emph{subgroup stratification} of $\Spec(R)$ as the stratification indexed by the poset of subgroups of $G$ whose $H$th stratum is the image of 
    \[ 
        \Spec(\Res_H^G R) \rightarrow \Spec(R) .
    \] 
\end{definition}

\begin{remark}
    There are various definitions of ``stratified space" in the literature. A common unifying theme of such definitions is that they form filtrations on a space with particularly nice properties. We opt to use ``stratification" over ``filtration" for two reasons. Firstly, when $H$ is a subgroup of $G$, one may refer to the $H$th stratum directly, whereas the alternatives in ``filtration" terminology would seem to be the relatively less readable ``$H$th filter" or ``$H$th filtered piece". 
    
    Secondly, in algebraic topology filtered objects tend to be filtered by infinite well--ordered sets, like $(\Z,\leq)$, rather than finite posets, whereas stratified spaces tend to have finitely many strata, for example in the theory of constructible sheaves. Moreover (co)filtered objects tend to have underlying $CW$-complexes, (like Postnikov/Whitehead towers, motivic and even filtrations, and slice towers), whereas stratified spaces tend to have more pathological algebro-geometric underlying point-set spaces.
\end{remark}

In all of our examples in \cref{sec:examples}, the strata happen to be closed subspaces.

\begin{question}
    Is the subgroup stratification always by closed subspaces?
\end{question}

Under somewhat strong assumptions on $G$ and on the Tambara functor $R$ we may establish further properties of the subgroup stratification.

\begin{definition}
    Let $G$ be a Dedekind group. We say that a Mackey, Green, or Tambara functor $M$ is \emph{idle} if, for each subgroup $H$ of $G$, the Weyl group $G/H$ acts trivially on $M(G/H)$.
\end{definition}

For example, the Burnside $G$-Tambara functor $\A_G$ is idle.

\begin{lemma}\label{lem:Nakaoka-spec-map-of-coind-unit-is-inj}
    Let $R$ be an idle Tambara functor and assume $G$ is Dedekind. The map 
    \[ 
        \Spec(\Res_H^G R) \cong \Spec(\Coind_H^G \Res_H^G R) \rightarrow \Spec(R)
    \] 
    induced by the adjunction unit $R \rightarrow \Coind_H^G \Res_H^G R$ is injective.
\end{lemma}

\begin{proof}
    Chasing through definitions, we see that $\p \mapsto \Res_H^G \p$ defines a set-theoretic left inverse, so long as we show $\Res_H^G \p$ is a prime ideal of $\Res_H^G R$. When checking whether or not an ideal is prime, it suffices to check merely the cases in which $I$ and $J$ are each generated by a single element. Suppose towards contradiction that there exist $x \in \Res_H^G R(H/L)$ and $y \in \Res_H^G R(H/L')$ such that $x, y \notin \Res_H^G \p$ and $(x) \cdot (y) \subset \Res_H^G \p$. 
    
    Let $x'$ denote the element of $R(G/L)$ corresponding to $x$ under the natural isomorphism $R(G/L) \cong \Res_H^G R(H/L)$, and similarly for $y'$. Since $R$ is idle, we have 
    \[ 
        \Res_H^G \left( (x') \cdot (y') \right) = \Res_H^G (x') \cdot \Res_H^G (y') = (x) \cdot (y) . 
    \]
    Finally, since $(x') \cdot (y')$ is generated as a Tambara ideal by elements in $R(G/L)$, where $L$ runs through subgroups of $H$, the statement 
    \[ 
        \Res_H^G \left( (x') \cdot (y') \right) = (x) \cdot (y) \subset \Res_H^G \p 
    \] 
    implies $(x') \cdot (y') \subset \p$, contradicting the fact that $\p$ is prime.
\end{proof}

It is a priori unclear that $\p \mapsto \Res_H^G \p$ respects primality of Tambara ideals (and in fact it does not in general: $\Res_e^G \p$ need only be a $G$-prime ideal of the $G$-ring $R(G/e)$). However the proof of \cref{lem:Nakaoka-spec-map-of-coind-unit-is-inj} shows that $\Res_H^G \p$ is prime, so long as $R$ is idle and $G$ is Dedekind. We thus immediately obtain the following.

\begin{corollary}\label{cor:idle+Dedekind-implies-stratification-is-injective}
    Assume $R$ is idle and $G$ is Dedekind. The assignment $\p \mapsto \Res_H^G \p$ determines an order-preserving poset morphism of Green ideals, which respects the property of being a Tambara ideal and respects primality. In particular, $\Spec(\Res_H^G R)$ is a retract (i.e. an idempotent subobject) of $\Spec(R)$ in the category of topological spaces.
\end{corollary}

\section{Domain-like Tambara functors}

As we have already built up enough technology to establish some results about domain-like Tambara functors, we digress on them here. This condition has not enjoyed much treatment in literature; in particular, while there are some results giving sufficient conditions for a Tambara functor to be domain-like, we give in \cref{thm:strata-usually-not-open-in-domain-like} a novel necessary condition, involving the subgroup stratification, for certain kinds of Tambara functors to be domain-like.

\begin{definition}[{\cite[Definition 4.28(1)]{Nak11a}}]\label{def:domain-like}
	A Tambara functor $R$ is called \emph{domain-like} if the zero ideal is prime.
\end{definition}

We have seen many times that the bottom level of a Tambara functor, viewed as a $G$-ring, carries much of the information of a Tambara functor. Following this idea, we set out to study the $G$-ring analogue of integral domain.

\begin{definition}
	An ideal $I$ of a $G$-ring $Q$ is $G$-invariant if $g(I) \subset I$ for all $g \in G$. We call $I$ \emph{$G$-prime} if furthermore $I \neq Q$ and whenever $x \cdot y \notin I$, there exists $g \in G$ such that either $gx \notin I$ or $gy \notin I$. If the zero ideal is $G$-prime, we call $Q$ domain-like.
\end{definition}

If $Q$ is a $G$-ring whose underlying ring is an integral domain, then $Q$ is domain-like. In other words, being domain-like is weaker than being a domain.

\begin{lemma}[{\cite[Lemma 4.3]{4DS}}]\label{lem:p-prime-implies-pGe-Gprime}
	Let $\p$ be a prime ideal of a $G$-Tambara functor $R$. Then $\p(G/e)$ is a $G$-prime ideal of $R(G/e)$.
\end{lemma}

In particular, if $R$ is a domain-like $G$-Tambara functor, then $R(G/e)$ is a domain-like $G$-ring. Now let $Q$ be a domain-like $G$-ring. The $G$-fixed points $Q^G$ are an integral domain, and $Q$ is an integral extension of $Q^G$ (since every element $\alpha$ of $Q$ is a root of the $G$-fixed monic polynomial $\prod_g (x-g\alpha)$).

The following is a slight strengthening of \cite[Proposition 4.31]{Nak11a}; note that Nakaoka also considers the condition on a Tambara functor in which all restrictions are injective, which in \cite[Definition 4.19]{Nak11a} is called the ``monomorphic restriction condition". In particular, we weaken the hypothesis from asking for $R(G/e)$ to be an integral domain to asking for $R(G/e)$ to be a domain-like $G$-ring.

\begin{lemma}\label{lem:domain-like-plus-MRC}
	Let $R$ be a Tambara functor such that $R(G/e)$ is a domain-like $G$-ring and all restrictions in $R$ are injective. Then $R$ is domain-like.
\end{lemma}

\begin{proof}
	Suppose that $I$ and $J$ are both not the zero ideal. Since restrictions are injective, both $I(G/e)$ and $J(G/e)$ are nonzero. Let $x \in I(G/e)$ and $y \in J(G/e)$ be nonzero. Since $R(G/e)$ is domain-like, there exists $g \in G$ such that $x \cdot gy$ is nonzero. Since $g \cdot J(G/e) = J(G/e)$, we see that $x \cdot gy$ is a nonzero element of $(I \cdot J)(G/e)$, so that $I \cdot J$ is nonzero. By contrapositive we see that the zero ideal is prime.
\end{proof}

\begin{proposition}\label{prop:domain-like-iff-coind-is}
	Let $R$ be either an $H$-Tambara functor or an $H$-ring. Then $\Coind_H^G R$ is domain-like if and only if $R$ is.
\end{proposition}

\begin{proof}
	We first prove the Tambara case. We observed in \cref{thm:coinduction-induces-bijection-on-Tambara-ideals} that the zero ideal is prime if and only if the coinduction of the zero ideal is prime. But the coinduction of the zero ideal is the zero ideal.
	
	Now let $Q$ be an $H$-ring. First, observe the natural isomorphism 
    \[ 
        \mathrm{FP}(\Coind_H^G Q) \cong \Coind_H^G \mathrm{FP}(Q) 
    \] 
    of $G$-Tambara functors. By \cref{lem:domain-like-plus-MRC,lem:p-prime-implies-pGe-Gprime}, $Q$ is domain-like if and only if $\mathrm{FP}(Q)$ is domain-like. By the Tambara functor case established in the previous paragraph, $\FP(Q)$ is domain-like if and only if $\Coind_H^G \mathrm{FP}(Q) \cong \mathrm{FP} (\Coind_H^G Q)$ is domain-like. But $\mathrm{FP} ( \Coind_H^G Q)$ is domain-like if and only if $\Coind_H^G Q$ is a domain-like $G$-ring.
\end{proof}

We give a sample application below.

\begin{example}\label{ex:some-domain-like-things}
	Let $\F$ be any field. First, form the $G$-ring $\F[G \cdot x]$, the free $\F$-algebra on $|G|$ generators which are permuted by $G$. Next, consider the $G$-ring $Q := \F[G \cdot x]/(g_1 x \cdot g_2 x)$, where $(g_1,g_2)$ runs through the complement of the diagonal in $G \times G$. By \cref{prop:domain-like-iff-coind-is}, $\Coind_e^G \F[t]$ is a domain-like Tambara functor, and we define a map $\mathrm{FP}(Q) \rightarrow \Coind_e^G \F[t]$ by sending $g x$ to the element which is $t$ in the $g^{\textrm{th}}$ factor (this is the adjoint of the ring map $Q \rightarrow \F[t]$ sending $e \cdot x$ to $t$ and $g x$ to zero for $g \neq e$).
	
	The ring map $Q \rightarrow \Coind_e^G \F[t]$ is visibly injective, hence the kernel of the Tambara functor map $\mathrm{FP}(Q) \rightarrow \Coind_e^G \F[t]$ is zero. Since the zero ideal in $\Coind_e^G \F[t]$ is prime, its preimage in $\mathrm{FP}(Q)$, which is also the zero ideal, is prime. Thus $\mathrm{FP}(Q)$ is domain-like. Additionally, $Q$ is a domain-like $G$-ring.
\end{example}

\begin{proposition}\label{prop:domain-like-iff-coind-from-clarified-domain-like}
	Let $R$ be a $G$-Tambara functor such that $R(G/e)$ is a Noetherian ring. $R$ is domain-like if and only if $R \cong \Coind_H^G R_H$ where $R_H$ is a clarified (\cref{def:clarified}) domain-like $H$-Tambara functor. It is possible that $H = G$, i.e. $R$ itself is already clarified.
\end{proposition}

\begin{proof}
	By our hypothesis, \cite[Theorem D]{Wis25a} applies, so that we may write 
    \[ 
        R \cong \prod_{H \subset G} \Coind_H^G R_H 
    \] 
    where each $R_H$ is clarified. In particular, we have 
    \[ 
        R(G/G) \cong \prod_{H \subset G} (\Coind_H^G R_H)(G/G) 
    \] 
    Suppose there exist $H \neq H'$ such that $R_H$ and $R_{H'}$ are both nonzero. Then $(\Coind_H^G R_H)(G/G)$ and $(\Coind_{H'}^G)(R_{H'})(G/G)$ are each nonzero. Let $d$ denote the unique idempotent of $R(G/G)$ such that 
    \[ 
        d \cdot R(G/G) = (\Coind_H^G R_H)(G/G)
    \] 
    and similarly for $d'$. Then $(d) \cdot (d')$ is the zero ideal, but neither $d$ nor $d'$ is zero. Consequently there is no such choice of $H$ and $H'$. The claim follows from \cref{prop:domain-like-iff-coind-is}.
\end{proof}

It seems like there are not simple conditions on even a clarified $G$-ring which are equivalent to being domain-like. For example, the $C_2$-ring $\mathbb{C}[x,y]/(xy)$, with $C_2$ acting by swapping $x$ and $y$ is a very simple domain-like $C_2$-ring (by \cref{ex:some-domain-like-things}). It is clarified, but not an integral domain. However, see \cite[Proposition 7.27]{4DS} for sufficient conditions when $G = C_p$ which seem to be necessary or almost necessary.

Knowing that a Tambara functor is domain-like is useful by the following special case of \cite[Theorem 4.9]{Nak11a}.

\begin{proposition}\label{prop:kernel-of-map-into-domain-is-prime}
	Let $R$ be a domain-like $G$-Tambara functor, and $f : S \rightarrow R$ a morphism of Tambara functors. Then the kernel of $f$ is a prime Tambara ideal.
\end{proposition}

\begin{proposition}\label{prop:inj-from-spec-to-spec}
	Let $R$ be any $G$-Tambara functor. There is an injection from the set of $G$-prime ideals of $R(G/e)$ to the set of prime Tambara ideals of $R$.
\end{proposition}

\begin{proof}
	Let $I$ be a $G$-prime ideal of $R(G/e)$, and let $J$ be the Tambara ideal which is the kernel of the Tambara functor morphism $R \rightarrow \mathrm{FP}(R/I)$ adjoint to the quotient $R(G/e) \rightarrow R(G/e)/I$. By \cref{lem:domain-like-plus-MRC,prop:kernel-of-map-into-domain-is-prime}, $J$ is prime. $I \mapsto J$ is the desired injection, because $I$ may be recovered from $J$ via the equality $J(G/e) = I$.
\end{proof}

We note that the image of the injection of \cref{prop:inj-from-spec-to-spec} is not in general the same as the $e$th stratum. Finally, we observe that the strata of domain-like Tambara functors fail to be open, in general.

\begin{theorem}\label{thm:strata-usually-not-open-in-domain-like}
    Let $G$ be an arbitrary group and $R$ any domain-like Tambara functor such that the restriction $R(G/G) \rightarrow R(G/H)$ is not injective. Then the $H$th stratum is not open.
\end{theorem}

\begin{proof}
    The closure of the zero ideal is all of $\Spec(R)$, so it suffices to show that the zero ideal is not contained in the $H$th stratum. The adjunction unit $R \rightarrow \Coind_H^G \Res_H^G R$ is given in level $G/G$ by the restriction $R(G/G) \rightarrow R(G/H)$. Thus no ideal in $\Coind_H^G \Res_H^G R$ can have preimage in $R$ equal to zero.
\end{proof}

\section{Explicit descriptions of some strata}\label{sec:examples}

\subsection{The Burnside Tambara functor}\label{subsec:Burnside_example}

We start by computing the stratification on the Nakaoka spectrum of the Burnside Tambara functor $\A_G$ when $G$ is Dedekind. Since $\A_G$ is idle, our strategy is to appeal to \cref{cor:idle+Dedekind-implies-stratification-is-injective}. We start by recalling the computation of the Nakaoka spectrum of $\A_G$, due in various cases to \cite{Nak14,CG23,4DSC}.

\begin{definition}
    Let $H \subset G$ and let $p$ denote either a prime integer or $0$. Define a ring homomorphism
    \[ 
        \phi^G_{H,p} : \A(G) \rightarrow \Z/p 
    \] 
    by $X \mapsto |X^H|$, for $X$ a $G$-set.
\end{definition}

\begin{definition}[{\cite[Definition 3.1]{CG23}}]
    Let $K$ and $H$ be subgroups of $G$, let $p$ denote either a prime integer or $0$, and define an ideal $\p_{K,p}$ of $\A_G$ by
    \[ 
        \p_{K,p}(G/H) := \bigcap_{\substack{I \subset H \\ I \preccurlyeq_G K}} \mathrm{Ker}(\phi_{H,p}^I)
    \] 
    where $I \preccurlyeq_G K$ means that the subgroup $I$ of $G$ is conjugate to a subgroup of $K$.
\end{definition}

\begin{theorem}[\cite{Nak14,CG23,4DSC}]
    Let $G$ be an arbitrary finite group. Every prime Tambara ideal of $\A_G$ is one of the $\p_{K,p}$. All inclusions and equalities between $\p_{K,p}$ are explicitly determined by the group structure of $G$.
\end{theorem}

Not all of the $\p_{K,p}$ are distinct. For example, by \cite[Figure 1]{CG23}, $\p_{C_2,2} = \p_{C_4,2}$ when $G = C_{12}$. In any case, we may now compute the strata of $\Spec(\A_G)$ when $G$ is Dedekind.

\begin{theorem}\label{thm:strata-of-Burnside-G-Dedekind}
    Let $G$ be a Dedekind group and $H \subset G$. The $H$th stratum of $\Spec(\A_G)$ is the subspace 
    \[ 
        \{ \p_{K,p} | K \subset H, p \in \Spec(\Z) \}
    \]
    and the map 
    \[ 
        \Spec(\A_H) \cong \Spec(\Coind_H^G \A_H) \rightarrow \Spec(\A_G) 
    \] 
    sends $\p_{K,p} \subset \A_H$ to $\p_{K,p} \subset \A_G$. The $H$th stratum of $\Spec(\A_G)$ is closed and not open.
\end{theorem}

\begin{proof}
    If $L \subset H \subset G$ and $K \subset G$, we have 
    \[ 
        \Res_H^G \p_{K,p}(H/L) = \p_{K,p}(G/L) = \bigcap_{I \subset L \cap K} \mathrm{Ker}(\phi_{L,p}^I) \subset \A(L) 
    \] 
    and
    \[ 
        \p_{K \cap H,p}(H/L) = \bigcap_{I \subset L \cap K \cap H} \mathrm{Ker}(\phi_{L,p}^I) \subset \A(L) .
    \] 
    Since $L \cap K \cap H = L \cap K$, we deduce
    \[ 
        \Res_H^G \p_{K,p} = \p_{K \cap H,p} .
    \]
    Thus by \cref{cor:idle+Dedekind-implies-stratification-is-injective} the map $\Spec(\A_H) \rightarrow \Spec(\A_G)$ has the claimed description on points. By \cref{thm:strata-usually-not-open-in-domain-like} the $H$th stratum is not open.

    It remains to show that the image of $\Spec(\A_H)$ in $\Spec(\A_G)$ is closed. First, since $\A_H$ is domain-like, the zero ideal $\p_{H,0}$ is prime. Second, recall the definition 
    \[ 
        V(\p_{H,0}) := \{ \p | \p_{H,0} \subset \p \} 
    \]
    Since $\Spec(\A_H) \rightarrow \Spec(\A_G)$ preserves inclusions of prime ideals (since $\Coind_H^G$ and preimage do), we see that the image of $\Spec(\A_H)$ in $\Spec(\A_G)$ is contained in the basic closed subset $V(\p_{H,0})$ of $\Spec(\A_G)$. To conclude the proof, we will show the opposite containment.

    Suppose $\p_{H,0} \subset \p_{K,p}$ for some possibly zero prime integer $p$. Note that $\p_{H,0}$ is the kernel of $\A_G \rightarrow \Coind_H^G \A_H$. Since $\A_H$ is, levelwise, a finitely-generated $\Z$-module, so is $\Coind_H^G \A_H$, so that $\A_G \rightarrow \Coind_H^G \A_H$ is levelwise integral. Therefore by the going-up theorem for Tambara functors \cite[Theorem 5.9]{4DS} we deduce that $\p_{K,p}$ is the preimage of some prime ideal $\Coind_H^G \q \subset \Coind_H^G \A_H$. It follows that $\q = \q_{L,q}$ for some possibly zero prime integer $q$ and $L \subset K$. By the first paragraph, we see that $\p_{K,p} = \q_{L,q}$.
\end{proof}

\begin{corollary}\label{cor:ex-of-etale-map-which-is-not-open}
    Let $G$ be an arbitrary finite group. There exists an \'{e}tale map of $G$-Tambara functors for which the induced map on Nakaoka spectra is not open. If $G$ is Dedekind, we may further construct an example which is closed.
\end{corollary}

\begin{proof}
    By \cite[Theorem A]{Wis25c} the map $\A_G \rightarrow \Coind_H^G \A_H$ is \'{e}tale, and the image on the induced map of Nakaoka spectra is precisely the $H$th stratum. If $H$ is a proper subgroup of $G$, the restriction $\A_G(G/G) \rightarrow \A_G(G/H)$ is not injective, so that the $H$th stratum is not open by \cref{thm:strata-usually-not-open-in-domain-like}.

    On the other hand, by the explicit description of the strata in \cref{thm:strata-of-Burnside-G-Dedekind}, we see that $\Spec(\Coind_H^G \A_H) \rightarrow \Spec(\A_G)$ is closed when $G$ is Dedekind.
\end{proof}

This is an extremely surprising result: classically, if $Q \rightarrow Q'$ is any \'{e}tale map of rings, the induced map $\Spec(Q') \rightarrow \Spec(Q)$ is open. It is natural to wonder whether or not this is a pathology of $\Coind_H^G$ specifically or of Tambara functors more generally. If $H$ is a subgroup of $G$ not equal to $G$, then $\Coind_H^G R$ is not clarified (unless $R$ is the zero $H$-Tambara functor). Thus in all of our examples of \'{e}tale maps of Tambara functors whose induced maps on Nakaoka spectra are not open, one of the two Tambara functors is not clarified.

\begin{question}
    Suppose $R \rightarrow S$ is an \'{e}tale map of clarified $G$-Tambara functors. Is 
    \[ 
        \Spec(S) \rightarrow \Spec(R) 
    \] 
    an open map?
\end{question}

\begin{remark}\label{rem:cf_4DSC_Theorem_4.19}
    By \cite[Proposition 4.19]{4DSC} there is a continuous bijection 
    \[ 
        \Phi : \Spec(\A_G(G/G)) \rightarrow \Spec(\A_G) 
    \] 
    from the Zariski spectrum of the Burnside ring $\A_G(G/G)$ to the Nakaoka spectrum of the Burnside Tambara functor; it is not a homeomorphism unless $G = e$. One might wonder whether the subgroup stratification is better behaved in $\Spec(\A_G(G/G))$ by taking the preimage, than in $\Spec(A_G)$. This is not the case, at least in the precise sense that for general finite groups $G$, the preimage of a general stratum along $\Phi$ remains non-open, by the following argument.

    Let $G$ be a group with no perfect subgroups other than $e$, and let $V$ denote the $H$th stratum of $\Spec(\A_G)$ for some subgroup $H$ of $G$ not equal to $G$. Since $V$ is closed and $\Phi$ is continuous, $\Phi^{-1} V$ is closed. Our assumption on $G$ implies $\Spec(\A_G(G/G))$ is connected by \cite[Theorem 2]{Sch79}. Our assumption on $H$ implies $V$ is not equal to $\Spec(\A_G)$, whence $\Phi^{-1} V$ is not equal to $\Spec(\A_G(G/G))$. Since $\Phi^{-1} V$ is a closed and proper nonempty subset of a connected space, it cannot be open.
\end{remark}

\subsection{The ghost}\label{subsec:ghost_example}

Recall from \cite[Definition 7.3]{4DS} that the \emph{ghost} of a $C_p$-Tambara functor $R$ is the $C_p$-Tambara functor $\ghost(R)$ given by the Lewis diagram
\[ \begin{tikzcd}
    R(C_p/e)^{C_p} \times R(C_p/C_p)/\Tr \arrow[d] \\
    R(C_p/e) \arrow[u, bend left = 45] \arrow[u, bend right = 45, swap]
\end{tikzcd} . \]
There is a canonical map $R \rightarrow \ghost(R)$ given by the identity in the $C_p/e$-level and the product of restriction and the quotient by the transfer ideal in the $C_p/C_p$-level. 

Now \cite[Theorem F]{4DS} says that every prime ideal of $\ghost(R)$ has the following description. Letting $(\p;\q)$ denote the ideal of $\ghost(R)$ specified by the Lewis diagram 
\[ \begin{tikzcd} 
    \p^{C_p} \times \q \arrow[d] \\
    \p \arrow[u, bend left = 45] \arrow[u, bend right = 45, swap]
\end{tikzcd} , \]
\cite[Theorem F]{4DS} says that every ideal of $\ghost(R)$ is either $(\p;R(C_p/C_p)/\Tr)$ for $\p$ a $C_p$-prime ideal of $R(C_p/e)$, or $(\Nm^{-1}\p;\p)$ for $\p$ a prime ideal of $R(C_p/C_p)/\Tr$.

\begin{proposition}\label{prop:stratum-of-ghost}
    Let $R$ be any nonzero $C_p$-Tambara functor and $\ghost(R)$ its ghost. 
    \begin{enumerate}
        \item A prime ideal $(\p;\q)$ of $\ghost(R)$ belongs to the $e$th stratum of $\Spec(\ghost(R))$ if and only if it is equal to $(\p;R(C_p/C_p)/\Tr)$ for $\p$ a prime ideal of $R(C_p/e)$.
        \item If $\Tr$ is not surjective then the $e$th stratum of $\Spec(\ghost(R))$ is not open.
        \item If $R(C_p/e)^{C_p}$ has finitely many minimal prime ideals (e.g. is Noetherian or an integral domain), the $e$th stratum of $\Spec(\ghost(R))$ is closed.
    \end{enumerate}
\end{proposition}

\begin{proof}
    The first statement follows from the observation that the map $\ghost(R) \rightarrow \Coind_e^{C_p} R(C_p/e)$ is given in level $C_p/C_p$ by the restriction 
    \[ 
        \ghost(R)(C_p/C_p) \rightarrow \ghost(R)(C_p/e) = R(C_p/e) . 
    \] 
    This restriction is the projection 
    \[ 
        \ghost(R)(C_p/C_p) = R(C_p/e)^{C_p} \times R(C_p/C_p)/\Tr \rightarrow R(C_p/e)^{C_p} . 
    \] 
    The preimage of an arbitrary prime ideal $\Coind_e^{C_p} \p$ of $\Coind_e^{C_p} R(C_p/e)$ is therefore $\bigcap_{g \in C_p} g \p$ in the $C_p/e$-level and $\p^{C_p} \times R(C_p/e)/\Tr$ in the $C_p/C_p$-level.

    For the second statement, let $\p$ be any prime ideal of $R(C_p/C_p)/\Tr$. Then $(\Nm^{-1}\p;\p)$ is not in the $e$th stratum, but the closure of $(\Nm^{-1}\p;\p)$ contains the point $(\Nm^{-1}\p;R(C_p/C_p)/\Tr)$ which is in the $e$th stratum. Thus the complement of the $e$th stratum is not closed.

    Finally, let $\p_i$ be the finite list of minimal prime ideals of $R(C_p/e)^{C_p}$. By the first paragraph, every prime ideal in the $e$th stratum has the form $(\q;R(C_p/C_p)/\Tr)$. Now for some $i$ we have $\p_i \subset \q$, whence $(\p_i;R(C_p/C_p)/\Tr) \subset (\q;R(C_p/C_p)/\Tr)$. Therefore the $e$th stratum is the finite union of the basic closed subsets $V((\p_i);R(C_p/C_p)/\Tr)$.
\end{proof}

\begin{corollary}\label{cor:ghost-strata-are-usually-not-open}
    Let $R$ be a $C_p$-Tambara functor such that $\Tr$ is not surjective. Then $\ghost(R) \rightarrow \Coind_e^{C_p} R(C_p/e)$ is \'{e}tale, but the induced map on Nakaoka spectra is not open.
\end{corollary}

\begin{proof}
    By \cite[Theorem A]{Wis25c} $\ghost(R) \rightarrow \Coind_e^{C_p} R(C_p/e)$ is \'{e}tale and by \cref{prop:stratum-of-ghost} the induced map on Nakaoka spectra is not open.
\end{proof}

Since the subgroup stratification is natural in Tambara functor morphisms and $\Spec(\ghost (R)) \rightarrow \Spec(R)$ is always surjective, \cref{prop:stratum-of-ghost} provides a method to determine the subgroup stratification on all $C_p$-Tambara functors.

\subsection{Fixed-point functors}

Let $Q$ be a $G$-ring. By \cite[Proposition 6.4]{4DS} the set of prime ideals of $\FP(Q)$ bijects with the set of prime ideals of $Q^G$, which is in bijection with the set of $G$-prime ideals of $Q$. Precisely, the proof of \cite[Proposition 6.4]{4DS} shows that this bijection sends a $G$-prime ideal of $Q$ to its intersection with $Q^G$, and its inverse sends a prime ideal of $Q^G$ to the unique $G$-prime ideal of $Q$ lying over it.

\begin{proposition}\label{prop:subgroup-strata-for-MRC}
    Let $Q$ be a $G$-ring. The $H$th stratum of $\Spec(\FP(Q))$ is the set of prime ideals $\p$ of $\FP(Q)$ such that the ideal $\p(G/H) \subset Q^H$ is prime.
\end{proposition}

\begin{proof}
    Let $\p$ be a prime ideal of $\FP(Q)$ such that $\p(G/H)$ is a prime ideal of $Q^H$. Then by \cite[Proposition 6.4]{4DS} there is a unique prime ideal $\q$ of $\FP(\Res_H^G Q) \cong \Res_H^G(\FP(Q))$ such that $\q(H/H) = \p(G/H)$. The preimage of $\Coind_H^G \q$ to $\FP(Q)$ is the unique prime ideal of $\FP(Q)$ given in level $G/G$ by the prime ideal $\q(H/H) \cap Q^G$. Since $\q(H/H) = \p(G/H)$ and $\p(G/G) = \p(G/H) \cap Q^G$, this preimage is precisely $\p$.
    
    Let $\p$ be a prime ideal of $\FP(Q)$ which is the preimage of a prime ideal $\Coind_H^G \q$ of $\Coind_H^G \Res_H^G \FP(Q)$. Since $\Res_H^G \FP(Q) \cong \FP(\Res_H^G(Q))$, \cite[Proposition 6.4]{4DS} implies $\q(H/H)$ is a prime ideal of $Q^H$. The $G/H$-level of the adjunction unit is the diagonal ring map
    \[ 
        \FP(Q)(G/H) \rightarrow (\Coind_H^G \Res_H^G \FP(Q))(G/H) \cong \FP(Q)(G/H)^{\times |G/H|} .
    \] 
    The preimage of $(\Coind_H^G \q)(G/H)$ along this map is simultaneously $\p(G/H)$ and the prime ideal $\q(H/H) \subset Q^H$.
\end{proof}

\begin{remark}
    The running hypothesis in \cite[Section 6]{4DS} that $R$ is a fixed-point Tambara functor may be weakened in some of their results to the assumption that all restrictions in $R$ are injective. Correspondingly we obtain a refinement of \cref{prop:subgroup-strata-for-MRC} to this case, with appropriate notational adjustments.
\end{remark}

\bibliographystyle{alpha}
\bibliography{ref}

\end{document}